\documentclass[axioms,article,accept,pdftex,moreauthors]{Definitions/mdpi}

\firstpage{1}
\pubvolume{1}
\issuenum{1}
\articlenumber{0}
\pubyear{2024}
\copyrightyear{2024}
\datereceived{ }
\daterevised{ } 
\dateaccepted{ }
\datepublished{ }
\hreflink{https://doi.org/} 

\Title{Some remarks on smooth mappings of Hilbert and Banach spaces and
	their local convexity property}

\TitleCitation{Some remarks on smooth mappings of Hilbert and Banach spaces and
	their local convexity property}

\Author{Yarema A. Prykarpatskyy $^{1, 2}$\orcidA{0000-0002-3033-7419},
Petro Ya. Pukach $^{3}$\orcidC{0000-0002-0359-5025},
Myroslava I. Vovk $^{4}$\orcidB{0000-0002-7818-7755},
	Michal Greguš ${^5}$\orcidD{0000-0002-8156-8962}
	 }

\AuthorNames{Yarema A. Prykarpatskyy, Petro Ya. Pukach, Myroslava I. Vovk,   Michal Greguš}

\address{%
$^{1}$ \quad Department of Applied Mathematics, University of Agriculture in Krakow, Al. Mickiewicza 21, 31-120,
Kraków, Poland; yarema.prykarpatskyy@urk.edu.pl\\
$^{2}$ \quad Institute of Mathematics of NAS of Ukraine, 3 Tereschenkivska st., Kyiv-4, 01024, Ukraine \\
$^{3}$ \quad Department of Computational Mathematics and Programming, Institute
of Applied Mathematics and Fundamental Sciences, Lviv Polytechnic National University, 12 Bandera Str., 79013 Lviv, Ukraine; petro.y.pukach@lpnu.ua \\
$^{4}$ \quad Department of Mathematics, Institute of Applied Mathematics and Fundamental Sciences, Lviv Polytechnic National University, 12 Bandera Str., 79013 Lviv, Ukraine; myroslava.i.vovk@lpnu.ua\\
$^{5}$ \quad Department of Information Systems, Faculty of Management, Comenius University in Bratislava, Odbojárov 10, 82005 Bratislava 25, Slovakia; michal.gregus@fm.uniba.sk

}

\corres{Correspondence: michal.gregus@fm.uniba.sk (M.G.)}

\abstract{We analyze smooth nonlinear mappings for Hilbert and Banach spaces that
carry small balls to convex sets, provided that the radius of the balls is
small enough. Being focused on the study of new and mild sufficient
conditions for a nonlinear mapping of Hilbert and Banach spaces to be
locally convex, we address a suitably reformulated local convexity problem
analyzed within \ the Leray-Schauder homotopy method approach for Hilbert
speces, and within the Lipscitz smoothness condition   both for Hilbert and
Banach spaces. Some of the results presented in the work prove to be interesting and
novel even for finite-dimensional problems. Open problems related to the
local convexity property for nonlinear mapping of Banach spaces are also
formulated.}


\begin{document}

\section{The local convexity mapping property}

\subsection{Introductory setting}

The importance of local convexity of nonlinear mappings of Hilbert and
Banach spaces in applied mathematics is well known \cite%
{AlTiFo,Aug,Frid-1,IoTi,Pr,Pryk,KZ,Pr1,Pr,BP,PBPP,SPS,Greg,Vovk}.
It is especially true of the theory of nonlinear differential-operator equations \cite%
{Che,Naz}, control theory and optimization theory\cite%
{Matv}, etc. Some of interesting and important for
applications local convexity properties, related with mappings of Hilbert
spaces, were first studied in \cite{Po}, and\ related with mappings of
Hilbert and Banach spaces were later generalized and studied both in \ \
\cite{BaSa,BoEmKo,Dyma,FaHa,Ivan,Ledy,Mord,Reis} and in \ \cite%
{BGHV,Ho,Go,Ge1,Ge2,Li,KZ,GH,MPVZ,Phel,PP}. It is worth reminding that a
nonlinear continuous mapping $f:X\rightarrow Y$ of Banach spaces $X$ and $Y$
is called to be \textit{locally convex}, if for any point $a\in X$ there
exists a ball $B_{\varepsilon }(a)\subset X$ of radius $\varepsilon >0,$
such that its image $f(B_{\varepsilon }(a))\subset Y$ is convex.

The property of local convexity holds for the special case of a
differentiable mapping $f:X\rightarrow Y$ of Hilbert spaces if the Frech\'{e}%
t derivative $f^{\prime }(x):X\rightarrow Y$ is Lipschitzian in a closed
ball $B_{r}(a)\subset X$ with radius $r>0$ which is centered at point $a\in
X $ and the linear mapping $f^{\prime }(a):X\rightarrow Y$, defined on the
whole space $X,$ closed and surjective. We need to mention here that this
notion differs from that introduced before in \ \ \cite{IoTi}. \ Concerning
a special case of a differentiable mapping $f:X\rightarrow Y$ of Hilbert
spaces, the property of local convexity, as it was \ first stated in \ \
\cite{Po}, holds and is based on the strong convexity of the ball $%
B_{r}(a)\subset X,$ \ if, in addition, the Frechet derivative $f_{x}^{\prime
}:X\rightarrow Y,x\in X,$ \ is Lipschitzian in a closed ball $%
B_{r}(a)\subset X$ of radius $r>0,$ centered at point $a\in X,$ and the
linear mapping $f_{a}^{\prime }:X\rightarrow Y$ is surjective. This
statement, is \ below proved making use of slightly different arguments from
those done before in \ \ \cite{BaBaPlPr,Po}, giving rise to an improved
estimation for the radius $\ $ of the ball $B_{\rho }(a)\subset X,\ $\ whose
image $f(B_{\rho }(a))\subset Y$ proves to be convex.\

More subtle techniques are required \cite{BGHV,GH,MPVZ} for the local
convexity problem for a nonlinear differentiable mapping $f:X\rightarrow Y$
of Banach spaces. Thus, its analysis is carried out only for the case of
Banach spaces with special properties. In particular, the locally convex
functions between Banach spaces are analyzed and the conditions guaranteeing
that a given function $f:U\rightarrow Y$ is locally convex are stated. As in
the Banach space case the local convexity property is based on much more
subtle properties both of a mapping and a Banach space under regard, special
attention is paid to the locally convex functions between Banach spaces. We
suitably reformulated local convexity problem for Banach spaces, taking into
account of the interplay between the modulus of convexity of a Banach space
and the modulus of smoothness of a function $f:X\rightarrow Y,$ before used
in the work \ \ \cite{BaBaPlPr} within a more general framework. \ It has
been noticed that some of the results presented in the preent work prove to
be interesting themselves and novel even for finite-dimensional problems.
Open problems related with the local convexity property for nonlinear
mapping of Banach spaces are also formulated.

Let $X,Y$ be Hilbert or Banach spaces and $U\subset X$ be an open set. \ \ A
mapping $f:U\rightarrow Y$ is said to be G\^{a}teaux differentiable \ \ \cite%
{Nir,Sch} \ at a point $a\in U,$ if there exists a continuous linear mapping
$f^{\prime }(a):X\rightarrow Y,$ such that for any $h\in X,$ there exists
the limit%
\begin{equation}
\lim_{t\rightarrow 0}\frac{f(a+th)-f(a)}{t}=f^{\prime }(a)h.  \label{B0a}
\end{equation}%
The mapping $f^{\prime }(a):X\rightarrow Y$ is called the \textit{G\^{a}%
teaux derivative} of the mapping $f:U\rightarrow Y$ $\ $at the point $a\in U.
$

A mapping $f:U\rightarrow Y$ $\ $is said to be Frechet differentiable at a
point $a\in U,$ if there exists a continuous linear mapping $f^{\prime
}(a):X\rightarrow Y,$ such that there exists the limit%
\begin{equation}
\lim_{||h||\rightarrow 0}\frac{||f(a+h)-f(a)-f^{\prime }(a)h||}{||h||}=0.
\label{B0b}
\end{equation}%
The mapping $f^{\prime }(a):X\rightarrow Y$ is called the \textit{Frechet
derivative} of the mapping $f:U\rightarrow Y$ $\ $at the point $a\in U.$ $\ $%
The $\ $mapping $f^{\prime }:X\rightarrow Y$ $\ $at the point $a\in U\ $is
\textit{Frechet differentiable}, it is also continuous and G\^{a}teaux
differentiable. A mapping $f:U\rightarrow Y$ is called \textit{continuously
Frechet differentiable, }if the mapping $X\ni a\rightarrow f^{\prime }(a)\in
L(X;Y)$ is \textit{continuous } with respect to the corresponding
topologies. Respectively, a mapping $f:U\rightarrow Y$ is called \textit{%
continuously G\^{a}teaux differentiable, }if the mapping $X\ni a\rightarrow
f^{\prime }(a)\in L(X;Y)$ is \textit{continuous } with respect to the
corresponding topologies. If a mapping $f:U\rightarrow Y$ is continuously G%
\^{a}teaux differentiable, then it is also continuously Frechet
differentiable, thus being of the\textit{\ }$C^{1}$- class. The latter means
that the simplest way to state the continuous Frechet differentiability of a
mapping $f:U\rightarrow Y$ is to state its continuous Frechet
differentiability. In general, a mapping $\ f:U\rightarrow Y$ is assumed to
be either continuously Frechet differentiable or locally Lipschitz. \ If a
continuously differentiable mapping $f:X\rightarrow Y$ is a bijection and
its inverse $f^{-1}:Y\rightarrow X$ is also continuously differentiable, it
is called a Frechet diffeomorphism. Taking into account the classical
Inverse Function Theorem \ \ \ \cite{Nir,Sch} \ , a continuously Frechet
differentiable mapping $f:X\rightarrow Y,$ such that for any $x\in X$ the
derivative $f^{\prime }(x):X\rightarrow Y$ is surjective, possesses the
inverse mapping $f^{\prime }(x)^{-1}:Y\rightarrow X,$ for which there exists
a positive constant $\alpha _{x}>0,$ such that $\ ||f^{\prime }(x)u||$ $\geq
\alpha _{x}||u||$ \ for any $u\in X,$ thus defining a local Frechet
diffeomorphism of Banach spaces $X$ and $Y.$ In particular, this means that
for each point $x\in X$ there exists an open set $U(x)\subset X,$ such that $%
f(U(x))\subset Y$ \ is open in $Y$ and $f|_{U(x)}:$ $U(x)\rightarrow f(U(x))$
is a local diffeomorphism. The general and often important for many
applications problem of describing conditions on a locally Frechet
diffeomorphic mapping $f:X\rightarrow Y$ to be globally diffeomorphic is
enough complicated, yet we are here interested only in the related and also
important for applications local convexity property of the constructed above
mapping $f|_{U(x)}:$ $U(x)\rightarrow f(U(x)),\ \ $ namely, in description
of conditions on the mapping $f:X\rightarrow Y,$ under which for any $x\in X$
and some small enough $\varepsilon >0$ \ there exists a ball $B_{\varepsilon
}(a_{x})\subset U(x),$ being a convex set centered at some point $a_{x}\in
U(x),$ and such that the image $f(B_{\varepsilon }(a_{x}))\subset Y\ $\ is
convex too. The local convexity problem sounds as follows:

\begin{Problem}
\label{Prob0} To \ construct at least sufficient conditions for a nonlinear
smooth mapping $f:X\rightarrow Y$ \ of a Hilbert\ space $X$ into a Hilbert
space $Y$ to be locally convex.
\end{Problem}

Let $(X,(\cdot |\cdot ))$ be a Hilbert space and a self-mapping  $%
f:X\rightarrow \ X$ is Frechet differentiable. The following proposition
demonstrates some typical  geometrical  conditions, which a priori can be
imposed on the self-mapping $f:X\rightarrow \ X$ \ \   to guarantee it to
satisfy  the local convexity property.

\begin{Proposition}
Let  $(X,(\cdot |\cdot ))$ be a Hilbert space and a   Frechet differentiable
mapping  $f:X\rightarrow X\ $ is such that the mapping $h:=f-cI:X\rightarrow
X,c\in \mathbb{R}\backslash \{0\},$ is proper, that is the preimage of a
compact set is compat too. Assume also, in addition  that if $U\subset X$ is
open and convex subset,   then the image $h(U)\subset X$ is simply connected
and for  any $\xi \in X$ $\ \sup_{u\in U}(\xi |u)<\infty ,$ then $\sup_{u\in
B_{\varepsilon }(a)}(\xi |h(u))<\infty .$ Then the mapping  $f:X\rightarrow X
$ is locally convex.
\end{Proposition}

\begin{proof}
To sketch a proof of this proposition, we take any $\theta \in \lbrack 0,1]$
and let $f_{\theta }:=\theta f+(1-\theta )I:X\rightarrow X$ and remark that
this mapping satisfies conditions assumed above too. Let $v_{0}\in f(U)$ and
$u_{0}\in f^{-1}(v_{0})$ and define the mapping $g_{\theta }(u):=f_{\theta
}(u)-f(u_{0}),$ $u\in U,$ satisfying the condition $g_{\theta }(\partial
U)\neq 0,$ where, obviously, $f_{\theta }(\partial U)=\partial f_{\theta
}(U),\theta \in \lbrack 0,1].$ If now easy to check that the Leray-Schauder
\ \ \cite{Nir,Sch,Pryk} \ topological degree $\deg (g_{\theta },U;0)$ is
well defined, the preimage $g_{\theta }^{-1}(0)\subset X$ is finite,
coinciding with the finite number of elements, as the Frechet derivative $%
g_{\theta }^{\prime }(u):X\rightarrow X$ is compact for any $u\in U.$ Taking
now into account that the mapping $[0,1]\ni \theta \rightarrow g_{\theta
}\in C^{1}(X)$ is a homotopy with $g_{0}=I-u_{0}$ and $g_{1}=f-v_{0},$ one
obtains easily that $\deg (f-v_{0},U;0)=1.$ Since the latter holds for every
$v_{0}\in f(U),$ it holds for all $\theta \in \lbrack 0,1],$ confirming that
the mapping $f|_{U}:U\rightarrow f(U)$ is a fiffeomorphism. To state that
the image $f(U)\subset X$ is convex, we assume that $v_{0}\in f(S),$ $\ $%
where $S\ \subset U$ is chosen in such a way that $f(S)\subset X$ is not a
strongly locally convex set. Let $\varepsilon >0$ be small enough positive
number, such that $\inf \{\lambda \in Spec(f^{\prime }(u):u\in
B_{\varepsilon }(u_{0})\}=\varepsilon _{0}>0.$ Then the mapping $%
h(u):=u_{0}+\varepsilon _{0}^{-1}(f(u)-f(u_{0})),u\in B_{\varepsilon
}(u_{0}),$ is invertible and satisfies the inequalities%
\begin{equation}
(\left( h^{-1}\right) ^{\prime }(v_{0})\alpha |\alpha )>\varepsilon
_{0}M_{0}^{-1}(\alpha |\alpha ),||\left( h^{-1}\right) ^{\prime
}(v_{0})\alpha ||\text{ }\leq ||\alpha ||  \label{B0d}
\end{equation}%
for any $\alpha \in X,$ where, by definiotion, $M_{0}:=\sup \{\lambda \in
Spec(f^{\prime }(u_{0}):u_{0}\in S\subset U\}.$ From \ inequality (\ref{B0d}%
) one derives that the vector $\ \left( h^{-1}\right) ^{\prime
}(v_{0})\alpha \in X$ lies in the interior of a strictly convex cone
symmetric around the axis passing through $v_{0}\in f(S)$ in the direction \
of a vector $\alpha \in X$ with respect the geometry induced on $X$ by the
inner product $(\cdot |\cdot ).$ Replacing the vector $\alpha \in X$ by $%
-\alpha \in X,$ we get that $-\left( h^{-1}\right) ^{\prime }(v_{0})\alpha
\in X$ lies in the interior of the strictly convex cone antipodal to that
described above. It means that for $\delta >0$ sufficiently small the points
$h^{-1}(p)\subset S$ and $p\in f(S)$ lie together in the interior of one of
these strictly convex cones, as well as the points $h^{-1}(q)\subset S$ and $%
q\in f(S)$ lie together in the interior of one of these strictly convex
antipodal cones. Based on these geometric properties followed by analogical
reasonings from \ \ \cite{Frid-2}, \ we finally state the convexity of the
image\ $f(U)\subset X.$
\end{proof}

\subsection{The local convexity property:  the Lipschitz smoothness analysis}

To enter now into the problem, we begin with a novel proof of the local
convexity result for a nonlinear smooth mapping $f:X\rightarrow Y\ $ of
Hilbert spaces, before formulated in \ \cite{BaBaPlPr,Po} and stated there
under different conditions. The following proposition holds.

\begin{Proposition}
\label{Pr_B1} Let $f:X\rightarrow Y$ \ \ be a nonlinear differentiable
mapping of Hilbert spaces whose Frechet derivative $\ f_{x}^{\prime
}:X\rightarrow Y,$ $x\in B_{r}(a),$ in a ball $B_{r}(a)\subset X$ centered
at point $a\in X,$ is Lipschitzian with a constant $L>0,$ the linear mapping
$f_{a}^{\prime }{}:X\rightarrow Y,a\in X,$ is surjective and the adjoint
mapping $f_{a}^{\prime ,\ast }{}:Y\rightarrow X$ satisfies the condition\ $%
||f_{a}^{\prime ,\ast }{}||$ $\geq \nu $ for some positive constant $\nu >0.$
Then for any $\varepsilon <\min \{r,\nu /(4L)$ the image $F_{\varepsilon
}(a):=$ $f(B_{\varepsilon }(a))\subset Y$ is convex.
\end{Proposition}

To prove Proposition \ref{Pr_B1}, it is useful to state the following simple
enough lemmas, based both on the Taylor expansion \cite%
{AlTiFo,FaHa,KZ,Nir,Sch} \ of the differentiable mapping $f:X\rightarrow Y$
at point $x_{0}\in B_{\varepsilon }(a)\subset B_{r}(a)$ and on the triangle
and parallelogram properties of the norm $||\cdot ||$ in a Hilbert space.

\begin{Lemma}
\label{Lm_B1a} Let a mapping $f_{x}^{\prime }:X\rightarrow Y,x\in X,$ \ be $%
L $-Lipschitzian in a ball $B_{\rho }(x_{0})\subset B_{\varepsilon }(a)$ of
radius $\rho >0,$ centered at point $x_{0}:=(x_{1}+x_{2})/2\in
B_{\varepsilon }(a)$ for arbitrarily chosen points $x_{1},x_{2}\in
B_{\varepsilon }(a).$ Then there exists such a positive constant $\mu >0$
that the\ norm $\ ||f_{x}^{\prime ,}{}^{\ast }(y)||$ $\geq \mu ||y||$ \ in
the ball $B_{\rho }(x_{0})\subset X$ \ for all $y\in Y,\ $ there holds the
estimation $||f(x_{0})-y_{0}||$ $\leq \rho \mu \ $\ for $%
y_{0}:=(y_{1}+y_{2})/2,$ $y_{1}:=f(x_{1}),y_{2}:=f(x_{2})$ and the equation $%
f(x)=y_{0}\ \ $possesses a solution $\bar{x}\in B_{\rho }(x_{0}),$ such that
$\ \ ||\bar{x}-x_{0}||$ $\leq \mu ^{-1}||f(x_{0})-y_{0}||.$
\end{Lemma}

\begin{proof}
Really, the \ following Taylor expansions at point $x_{0}\in B_{\varepsilon
}(a)$ hold:

\begin{eqnarray}
y_{1} &=&f(x_{1})=f(x_{0})+f_{x_{0}}^{\prime }(x_{1}-x_{0})+\epsilon _{1},
\label{B1} \\
&&  \notag \\
y_{2} &=&f(x_{2})=f(x_{0})+f_{x_{0}}^{\prime }(x_{2}-x_{0})+\epsilon _{2},
\notag
\end{eqnarray}%
where $||\epsilon _{j}||$ $\leq \frac{L}{2}||x_{j}-x_{0}||^{2}=\frac{L}{8}%
||x_{1}-x_{2}||^{2},j=\overline{1,2},$ as the mapping $f_{x}^{\prime
}:X\rightarrow Y,$ $x\in X,$\ is $L$-Lipschitzian. From \ (\ref{B1}) one
obtains easily that%
\begin{equation}
y_{0}=f(x_{0})+\epsilon _{0},  \label{B1a}
\end{equation}%
where, evidently, $||\epsilon _{0}||$ $\leq (||\epsilon _{1}||+||\epsilon
_{2}||)/2\leq \frac{L}{8}||x_{1}-x_{2}||^{2}.$ Moreover, owing to the
Lipschitzian property of the Frechet derivative $f^{\prime }(x):X\rightarrow
Y,$ \ one can obtain the following inequality:
\begin{equation}
\begin{array}{c}
||f_{x}^{\prime ,\ast }(y)||=||f_{x}^{\prime ,\ast }(y)-f_{a}^{\prime ,\ast
}(y)+f_{a}^{\prime ,\ast }(y)||\text{ }\geq  \\
\\
\geq ||f_{a}^{\prime ,\ast }(y)||-||f_{x}^{\prime ,\ast }(y)-f_{a}^{\prime
,\ast }(y)||\text{ }\geq  \\
\\
\geq \nu ||y||-L||x-a||\cdot ||y||\text{ }\geq (\nu -L\varepsilon
)||y||:=\mu _{0}||y||%
\end{array}
\label{B1b}
\end{equation}%
for $\mu _{0}=(\nu -L\varepsilon )>0,$ as the norm $||x-a||$ $\leq
\varepsilon .$ This, in particular, means that the adjoint mapping $%
f_{x}^{\prime ,\ast }:Y\rightarrow X\ $\ at $x\in B_{\varepsilon }(a)$ \ is
invertible, defined on the whole Hilbert space $Y$ and the norm of its
inverse mapping $(f_{x}^{\prime ,\ast })^{-1}:X\rightarrow Y$ \ is bounded
on the ball $B_{\varepsilon }(a)\subset X$ by the value $1/\mu _{0}.$ First
observe that for $\rho :=\frac{\mu _{0}}{8\varepsilon \mu }%
||x_{1}-x_{2}||^{2}$ and $\mu :=\nu -2L\varepsilon >0$ the following
inequality
\begin{equation}
\begin{array}{c}
||f(x_{0})-y_{0}||=||\epsilon _{0}||\text{ }\leq \frac{L}{8}%
||x_{1}-x_{2}||^{2}= \\
\\
=\rho \mu \mu _{0}^{-1}L\varepsilon \leq \rho \mu \mu _{0}^{-1}(\nu
-L\varepsilon )-\rho \mu \mu _{0}^{-1}(\nu -2L\varepsilon )\leq  \\
\\
\leq \rho \mu \mu _{0}^{-1}(\nu -\ L\varepsilon )=\rho \mu ,%
\end{array}
\label{B1c}
\end{equation}%
based on \ expression (\ref{B1a}), \ holds. \ \ Denote now by $\bar{x}$ $\in
B_{\varepsilon }(a)$ an arbitrary point satisfying the condition $y_{0}=f(%
\bar{x}),$ whose existence is guaranteed by the standard implicit function
theorem \cite{Sch,Nir}, and denote $\bar{y}:=\left(
\int_{0}^{1}f_{x(t)}^{\prime ,\ast }dt\right) ^{-1}(\bar{x}-x_{0})\ $ $\in Y,
$ \ where the linear mapping $\left( \int_{0}^{1}f_{x(t)}^{\prime ,\ast
}dt\right) ^{-1}\ :X\rightarrow Y$ \ is bouned and \ determined owing to the
homotopy equality
\begin{equation}
f(\bar{x})-f(x_{0})=\int_{0}^{1}f_{x(t)}^{\prime }(\bar{x}-x_{0})dt:=\left(
\int_{0}^{1}f_{x(t)}^{\prime }dt\right) (\bar{x}-x_{0}),  \label{B1d}
\end{equation}%
which holds owing to the continuation $x(t):=x_{0}+t(\bar{x}-x_{0})\in
B_{2\varepsilon }(x_{0})\ $ for $t\in \lbrack 0,1].$ Moreover, the following
estimation
\begin{equation}
\left\Vert \left( \int_{0}^{1}f_{x(t)}^{\prime ,\ast }dt\right)
^{-1}\right\Vert \leq \mu ^{-1}  \label{B1e}
\end{equation}%
holds. Really, for any $y\in Y$ \
\begin{equation}
\begin{array}{c}
\begin{array}{c}
\left\Vert \left( \int_{0}^{1}f_{x(t)}^{\prime ,\ast }dt\right)
(y)\right\Vert =\left\Vert \left( \int_{0}^{1}f_{x_{0}}^{\prime ,\ast
}dt\right) (y)+\left( \int_{0}^{1}\left( f_{x(t)}^{\prime ,\ast
}-f_{x_{0}}^{\prime ,\ast }dt\right) dt\right) (y)\right\Vert \geq  \\
\\
\geq \left\Vert \left( \int_{0}^{1}f_{x_{0}}^{\prime ,\ast }dt\right)
(y)\right\Vert -\ \left\Vert \left( \int_{0}^{1}\left( f_{x(t)}^{\prime
,\ast }-f_{x_{0}}^{\prime ,\ast }dt\right) dt\right) (y)\right\Vert \geq  \\
\\
\geq \mu _{0}||y||-\frac{L}{2}||\bar{x}-x_{0}||\text{ }||y||\text{ }\geq
\left( \mu _{0}-L\varepsilon \right) ||y||\text{ }=\left( \nu -2L\varepsilon
\right) ||y||=\mu ||y||,%
\end{array}%
\end{array}
\label{B1ee}
\end{equation}%
as the norm $||\bar{x}-x_{0}||=||(\bar{x}-a)+(a-x_{0})||$ $\leq $ $||(\bar{x}%
-a)||+||(a-x_{0})||$ $\leq 2\varepsilon .$

$\ $ \ \ \textit{Remark}. The integral expressions, considered \ above with
respect to the parameter $t\in \lbrack 0,1],$ are well defined in a Hilbert,
or in general, in a Banach space $Y,$ as the\ related mapping $f_{x(\circ
)}^{\prime }(\bar{x}-x_{0}):[0,1]\rightarrow Y,$ being continuous and of
bounded variation, is \ a priori Riemann-Birkhoff type integrable \ \ \cite%
{BaMu,BaPo,Birk,CaCr,DiUh,KaSh,Jake}.

Denote by \ $(\cdot |\cdot )$ the scalar product both on the Hilbert space $%
X $ and the Hilbert space $Y.$ Then one easily obtains that
\begin{equation*}
\begin{array}{c}
||\bar{x}-x_{0}||^{2}=|(x-x_{0}|\left( \int_{0}^{1}f_{x(t)}^{\prime ,\ast
}dt\right) (\bar{y}))|=|(\int_{0}^{1}f_{x(t)}^{\prime }(\bar{x}-x_{0})|\bar{y%
})|\text{ }= \\
\\
=|(f(\bar{x})-f(x_{0})|\bar{y})|\text{ }=|(y_{0}-f(x_{0})|\left(
\int_{0}^{1}f_{x(t)}^{\prime ,\ast }dt\right) ^{-1}(\bar{x}-x_{0}))|\text{ }%
\leq \\
\\
\leq ||y_{0}-f(x_{0})||\text{ }||\left( \int_{0}^{1}f_{x(t)}^{\prime ,\ast
}dt\right) ^{-1}||\text{ }||\bar{x}-x_{0}||\text{ }\leq \\
\\
\leq ||y_{0}-f(x_{0})||\text{ }||\bar{x}-x_{0}||/\mu ,%
\end{array}%
\end{equation*}%
yielding the searched for inequality
\begin{equation}
||\bar{x}-x_{0}||\text{ }\leq \ \mu ^{-1}||f(x_{0})-y_{0}||,  \label{B1f}
\end{equation}%
and proving the Lemma.

\begin{Lemma}
\label{Lm_B1b}For arbitrarily chosen points $x_{1},x_{2}\in B_{\varepsilon
}(a)$ the whole ball $B_{\rho }(x_{0})$ of radius $\rho =\frac{\mu _{0}}{%
8\varepsilon \mu }||x_{1}-x_{2}||^{2}\leq \varepsilon ,$ centered at point $%
x_{0}:=(x_{1}+x_{2})/2\in B_{\varepsilon }(a),$ $\ $ belongs to the ball $%
B_{\varepsilon }(a).$
\end{Lemma}
\end{proof}

\begin{proof}
Consider for this the following triangle inequality and the related
parallelogram identity on the Hilbert space $X$ for any point $x\in B_{\rho
}(x_{0}):$%
\begin{equation}
\begin{array}{c}
||x-a||=||(x-x_{0})+(x_{0}-a)||\text{ }\leq ||x-x_{0}||\text{ }+||x_{0}-a||%
\text{ }= \\
\\
=||x-x_{0}||\text{ }+\text{ }||(x_{1}-a)/2+(x_{2}-a)/2||\text{ }= \\
\\
=||x-x_{0}||\text{ }%
+[(||x_{1}-a||^{2}/2+||x_{2}-a||^{2}/2)-||x_{1}-x_{2}||^{2}/4]^{1/2}\leq \\
\\
\leq \rho +(\varepsilon ^{2}-||x_{1}-x_{2}||^{2}/4)^{1/2}.%
\end{array}
\label{B1g}
\end{equation}%
For the righthand side of \ (\ref{B1g}) to \ be equal or less of $%
\varepsilon >0,$ \ it is enough to take such a positive number $\rho \leq
\varepsilon $ that

\begin{equation}
\rho +(\varepsilon ^{2}-||x_{1}-x_{2}||^{2}/4)^{1/2}\leq \varepsilon .
\label{B1h}
\end{equation}%
This means that the following inequality should be satisfied:%
\begin{equation}
\rho ^{2}\geq 2\varepsilon \rho -||x_{1}-x_{2}||^{2}/4.  \label{B1i}
\end{equation}%
The choice $\rho =\frac{\mu _{0}}{8\varepsilon \mu }||x_{1}-x_{2}||^{2}\
\geq 2L\varepsilon ^{2}/(\nu -\varepsilon L)$ \textit{a priori }satisfies
the above condition \ (\ref{B1i}) if $\varepsilon \leq \nu /(3L),$ \ thereby
proving the Lemma.
\end{proof}

\begin{proof}
(\textit{Proof of Proposition \ref{Pr_B1}).} Now, \ based on Lemmas \ \ref%
{Lm_B1a} and \ \ref{Lm_B1b}, it is easy to observe from \ (\ref{B1c}) and \
\ (\ref{B1f}) that a point $\bar{x}\in B_{\varepsilon }(a),$ satisfying the
equation $\ y_{0}=f(\bar{x}),$ belongs to the ball $\ B_{\rho
}(x_{0})\subset X:$
\begin{equation}
||\bar{x}-x_{0}||\text{ }\leq ||y_{0}-f(x_{0})||/\mu \leq \rho \mu /\mu
=\rho ,  \label{B2}
\end{equation}%
thereby proving our Proposition \ref{Pr_B1} and solving Problem (\ref{Prob0}%
).
\end{proof}

It is worth to mention here that our local convexity proof for a nonlinear
smooth mapping $f:X\rightarrow Y$ of Hilbert spaces is slightly different
from that presented before in \ \ \cite{Po} and gives a bit improved
estimation of the radius $\ $ of the ball $B_{\rho }(x_{0})\subset X,$ whose
image $f(B_{\rho }(x_{0}))\subset Y$ proves to be convex.

\section{The local convexity mapping property: Banach spaces case}

\subsection{Banach spaces of the modulus of convexity of degree 2}

Let now $X,Y$ be Banach spaces. A mapping $f:U\rightarrow Y$ defined on an
open subset $U\subset X$ we will call \emph{locally convex} if for each
point $x\in U$ and its neighborhood $O_{x}\subset U$ \ there is a convex
open neighborhood $U_{x}\subset O_{x}$ with convex image $f(U_{x})\subset Y.$
In this part we address the following problem.

\begin{Problem}
\label{Prob1} Find at least sufficient conditions guaranteeing that a given
function $f:U\rightarrow Y$ is locally convex.
\end{Problem}

We start from the modulus convexity definition following \cite{FaHa} \ (se
also \cite{MPVZ}).

\begin{Definition}
Let $(X,||\cdot ||)$ be a Banach space and $B_{X}:=\{x\in X:||x||$ $\leq
1\}\ $ be a unit ball. \ For every $\varepsilon $ $\in (0,2]$ we define the
modulus \ of convexity (or rotundity) of $||\cdot ||$ \ by \
\begin{equation}
\delta _{X} (\varepsilon)=\inf_{x,y\ \in B_{X}}\left\{ 1-\left\Vert \frac{x+y%
}{2}\right\Vert :\left\vert \left\vert x-y\right\vert \right\vert \geq
\varepsilon \right\} .  \label{3.1}
\end{equation}%
The norm $||\cdot ||$ \ is called uniformly convex (or uniformly rotund), if
$\delta _{X}\{\varepsilon )$ $>0$ for all $\varepsilon $ $\in (0,2].$ The
space $(X,||\cdot ||)$ is then called a uniformly convex space. \ Note also
that $\delta _{X}\{\varepsilon )=\inf_{Y\subset X}\{\delta _{Y}\{\varepsilon
):\dim Y=2\}.$
\end{Definition}

It is easy to observe that $\delta _{X}(\varepsilon )\leq \varepsilon /2$ \
for all $\varepsilon $ $\in (0,2].$ \ The definition above can be
equivalently reformulated owing to the following\ \ \ \cite{FaHa} \ lemma.

\begin{Lemma}
\ \label{Lm_1} $(X,||\cdot ||)$ be a Banach space, $S_{X}:=\partial B_{X}$
be the boundary of the unit ball $B_{X}\subset X$ and let $\delta
_{X}\{\varepsilon ),$ $\varepsilon $ $\in (0,2],$ be the modulus of
convexity \ of $||\cdot ||.$ \ \ Then
\begin{equation}
\delta _{X} (\varepsilon )=\inf_{x,y\ \in S_{X}}\left\{ 1-\left\Vert \frac{%
x+y}{2}\right\Vert :\left\vert \left\vert x-y\right\vert \right\vert
=\varepsilon \right\} .  \label{3.2}
\end{equation}
\end{Lemma}

Remark also, that by \cite{MPVZ}, each Banach space with norm, having
modulus of convexity of power type 2, is superreflexive. In addition, based
on \ the definition \ (\ref{3.1}) one can be derived \ \ \ \cite{FaHa} \ the
following useful lemma.

\begin{Lemma}
\ \label{Lm_2} The norm of a Banach space $X$ has modulus of convexity of
power $p>1,$ if and only if there is a positive constant $C>0$ such that
\begin{equation}
\left\Vert \frac{x+y}{2}\right\Vert \leq 1-C\Vert x-y\Vert ^{p}  \label{3.2a}
\end{equation}%
for any points $x,y\in X$ with $\max \{\Vert x\Vert ,\Vert y\Vert \}\leq 1.$
\end{Lemma}

\subsection{The Banach space case: main result}

To answer the Problem~above, we need to recall some notions related to the
differentiability and the Lipschitz property.

Let $X,Y$ \ be Banach spaces and $U\subset X$ be an open subset in $X.$ A
function $f:U\rightarrow Y$ is called

\begin{itemize}
\item \emph{differentiable} at a point $x_{0}\in U$ if there is a linear
continuous operator $f_{x_{0}}^{\prime }:X\rightarrow Y$ (called the Frechet
\emph{derivative} of $f$ at $x_{0}\in U$) such that
\begin{equation*}
\lim_{x\rightarrow x_{0}}\frac{\Vert f(x)-(f(x_{0})+f_{x_{0}}^{\prime
}(x-x_{0}))\Vert }{\Vert x-x_{0}\Vert }=0;
\end{equation*}

\item \emph{Lip-differentiable} at a point $x_{0}\in U,$ if there is a
neighborhood $W\subset U$ of $x_{0}\in U,$ such that $f:U\rightarrow Y$ is
differentiable at each point $x\in W$ and
\begin{equation*}
\sup_{x,y\in W,x\neq y}\frac{\Vert f(y)-(f(x)+f_{x}^{\prime }(y-x))\Vert }{%
\Vert y-x\Vert ^{2}}<\infty .
\end{equation*}

\item \emph{locally Lipschitz} at a point $x_{0}\in U,$ if there is a
neighborhood $W\subset U$ of $x_{0}\in U,$ such that
\begin{equation*}
\sup_{x,y\in W,x\neq y}\frac{\Vert f(x)-f(y)\Vert }{\Vert x-y\Vert }<\infty .
\end{equation*}
\end{itemize}

\begin{Theorem}
\ \ \label{Tm_3} Let $X$ be a Banach spaces whose norm has modulus of
convexity of power type 2. A homeomorphism $f:U\rightarrow V$ between two
open subsets $U,V\subset X$ is locally convex if

\begin{itemize}
\item the function $f:U\to V$ is Lip-differentiable at each point $x_0\in U$
and

\item the function $f^{-1}:V\rightarrow U$ is locally Lipschitz at each
point $y_{0}\in V.$
\end{itemize}
\end{Theorem}

\begin{proof}
Fix any point $x_{0}\in U.$ \ Given a neighborhood $O(x_{0})\subset U$ of $%
x_{0}\in U$ we should construct a convex neighborhood $U(x_{0})\subset
O(x_{0})$ with convex image $f(U(x_{0})).$ We lose no generality assuming
that $y_{0}=f(x_{0})=0.$

Using the Lip-differentiability of $f:U\rightarrow V$ at $x_{0}\in U,$ find
a neighborhood $W\subset O(x_{0})$ of $x_{0}\in U$ and a real number $L$
such that
\begin{equation*}
\Vert f(y)-f(x)-f_{x}^{\prime }(y-x)\Vert \leq L\Vert y-x\Vert ^{2}
\end{equation*}%
for all points $x,y\in W.$ Moreover, since the homeomorphism $%
f^{-1}:V\rightarrow U$ is locally Lipschitz at the point $y_{0}=f(x_{0}),$
we can assume that
\begin{equation*}
\Vert x-y\Vert =\Vert f^{-1}(f(x)-f^{-1}(f(y))\Vert \leq L\Vert
f(y)-f(x)\Vert
\end{equation*}%
for any points $x,y\in W.$ \ \ By our assumption, the norm of the Banach
space $X$ \ has modulus of convexity of power type $2.$ \ Then, owing to the
result (\ref{3.2a}) of Lemma \ref{Lm_2}, there is a positive constant $C<1$
such that
\begin{equation}
\frac{1}{2}\ ||{x+y||}\ \leq 1-C\Vert x-y\Vert ^{2}  \label{3.3}
\end{equation}%
for any points $x,y\in B_{X}$ $.$ Take any $\varepsilon >0,$ such that

\begin{itemize}
\item $\frac{1}{4}L^{2}\leq \frac{C}{\varepsilon };$

\item $B_{\varepsilon }(z_{0})=\{z\in X:\Vert z-z_{0}\Vert \leq \varepsilon
\}\subset f(W);$
\end{itemize}

The choice of $\varepsilon >0$ guarantees that the preimage $A_{\varepsilon
}=f^{-1}(B_{\varepsilon }(z_{0}))$ of the $\varepsilon $-ball $%
B_{\varepsilon }(z_{0})=\{z\in V:\Vert z_{0}-z\Vert $ $\leq \varepsilon
\}=\{z\in X:\Vert z\Vert $ $\leq \varepsilon \},$ centered at $z_{0}=0,$
lies in the neighborhood $W$ on the point $x_{0}\in U.$ Now the proof of the
theorem will be complete as soon as we check that the closed neighborhood $%
A_{\varepsilon }=f^{-1}(B_{\varepsilon }(z_{0}))\subset W$ of $x_{0}\in U$
is convex. It suffices to check that for any points $x,y\in A_{\varepsilon }$
the point $\bar{x}=(x+y)/2$ belongs to $A_{\varepsilon }\subset W,$ which
happens if and only if its image $f(\bar{x})\in Y$ \ belongs to the ball $%
B_{\varepsilon }(z_{0})\subset Y.$ The choice of $f:U\rightarrow V$
guarantees that
\begin{equation*}
\Vert f(x)-f(\bar{x})-f_{z}^{\prime }(x-\bar{x})\Vert \leq L\Vert x-\bar{x}%
\Vert ^{2}
\end{equation*}%
and
\begin{equation*}
\Vert f(y)-f(\bar{x})-f_{z}^{\prime }(y-\bar{x})\Vert \leq L\Vert y-\bar{x}%
\Vert ^{2}.
\end{equation*}%
Adding these inequalities and taking into account that $x-\bar{x}=-(y-\bar{x}%
),$ we get
\begin{equation*}
\begin{array}{c}
\Vert (f(x)+f(y))-2f(\bar{x})\Vert \leq 2L\Vert y-\bar{x}\Vert ^{2}= \\
\\
=\frac{1}{2}L\Vert y-x\Vert ^{2}\leq \frac{1}{2}L^{2}\Vert f(y)-f(x)\Vert
^{2}%
\end{array}%
\end{equation*}%
and hence
\begin{equation*}
\
\begin{array}{c}
\left\Vert \frac{f(x)+f(y)}{2}-f(\bar{x})\right\Vert \ \leq \frac{1}{4}%
L^{2}\Vert f(x)-f(y)\Vert ^{2}\leq \\
\\
\leq \frac{C}{\varepsilon }\Vert f(x)-f(y)\Vert ^{2}.%
\end{array}%
\end{equation*}%
Since $f(x),f(y)\in B_{\varepsilon }(z_{0}),$ we get $\left\Vert \frac{%
f(x)+f(y)}{2}\right\Vert \leq \varepsilon -\frac{C}{\varepsilon }\Vert
f(y)-f(x)\Vert ^{2}$ and
\begin{eqnarray*}
\Vert f(z)\Vert &\leq &\left\Vert \frac{f(x)+f(y)}{2}\right\Vert +\left\Vert
\frac{f(x)+f(y)}{2}-f(z)\right\Vert \leq \\
&& \\
&\leq &\varepsilon -\frac{C}{\varepsilon }\Vert f(y)-f(x)\Vert ^{2}+\frac{C}{%
\varepsilon }\Vert f(y)-f(x)\Vert ^{2}=\varepsilon ,
\end{eqnarray*}%
which means that $f(z)\in B_{\varepsilon }(y_{0})$ and $z\in A_{\varepsilon
}\subset W.$
\end{proof}

\textit{Remark. }As follows from the proof of Theorem \ref{Tm_3}, \ when the
spaces $X$ and $Y$ are Hilbert ones, \ the obtained above result \ reduces
to that of Proposition \ref{Pr_B1}, yet in its slightly weakened form.

\subsection{The locally convex functions between Banach spaces}

As above, let $X,Y$ be Banach spaces and \ a function $f:U\rightarrow Y$ be
defined on an open subset $U\subset X,$ which is called \emph{locally convex,%
} if each point $x\in U$ has a neighborhood base consisting of open convex
subsets $U_{x}\subset X$ with convex images $f(U_{x})\subset Y.$ \ In this
part of paper we address the local convexity problem problem (\ref{Prob1})
formulated before. \ The answer to this problem will be given in terms of
the interplay between the modulus of smoothness of the function $\ \
f:U\rightarrow Y$ and the modulus of convexity of the Banach space $Y.$

Any Hilbert space $E$ of dimension $\dim (E)>1$ has modulus of convexity $%
\delta _{E}(t)=1-\sqrt{1-t^{2}/4}\leq \frac{1}{8}t^{2}.$ \ By \cite{GH}, $%
\delta _{X}(t)\leq \delta _{E}(t)\leq \frac{1}{8}t^{2}$ for each Banach
space $X.$ We shall say that the Banach space $X$ has \emph{modulus of
convexity of degree $p$} if there is a constant $L>0$ such that $\delta
_{X}(t)\geq Lt^{p}$ for all $t\in \lbrack 0,2]\ $that \ follows from the
inequlities $Lt^{p}\leq \delta _{X}(t)\leq \frac{1}{8}t^{2}$ that $p\geq 2.$
So, the Hilbert spaces have modulus of convexity of degree $2.$ Many
examples of Banach spaces with modulus of convexity of degree $2$ can be
found in \cite{HMZ}. In particular, so is the Banach space $%
X=(\sum\limits_{n=1}^{\infty }l_{4}(n))_{l_{2}},$ which is not isomorphic to
a Hilbert space. \medskip

Next, we recall \cite{GH,HMZ} \ the definition of the moduli of smoothness $%
\omega _{n}(f;t),t\geq 0,$ of a function $f:U\rightarrow Y$ defined on a
subset $U\subset X$ of a Banach space $X.$ By definition,
\begin{equation}
\begin{array}{c}
\omega _{n}(f;t)=\sup \Big\{\big\|\sum_{k=0}^{n}(-1)^{n-k}\binom{n}{k}f(x+
\\
\\
+(\tfrac{n}{2}-k)h)\big\|:\Vert h\Vert \leq t,\;[x-nh/2,x+nh/2]\subset U%
\Big\}.%
\end{array}
\label{3.4}
\end{equation}%
In particular,
\begin{equation*}
\omega _{1}(f;t)=\sup \{\Vert f(x)-f(y)\Vert :\Vert x-y\Vert \leq
t,\;[x,y]\subset U\}
\end{equation*}%
and
\begin{equation*}
\omega _{2}(f;t)=\sup \{\Vert f(x+h)-2f(x)+f(x-h)\Vert :\Vert h\Vert \leq
t,\;[x-h,x+h]\in U\},
\end{equation*}%
where $[x,y]=\{sx+(1-s)y:s\in \lbrack 0,1]\}$ stands for the segment
connecting the points $x,y\in X.$ Moreover, it is true \cite{GH} that $%
\omega _{n}(f;\frac{1}{m}t)\geq \frac{1}{m^{n}}\omega _{n}(f;t)$ for each $%
m\in \mathbb{N},t\geq 0.$ Below we will formulate the following definition.

\begin{Definition}
We shall say that a function $f:U\rightarrow Y,U\subset X,$ is \
\end{Definition}

\begin{itemize}
\item \emph{Lipschitz} if there is a constant $L$ such that $\omega
_{1}(f;t)\leq Lt$ for all $t\geq 0;$

\item \emph{second order Lipschitz} if there is a constant $L$ such that $%
\omega _{2}(f;t)\leq Lt^{2}$ for all $t\geq 0;$

\item \emph{locally} (\emph{second order}) \emph{Lipschitz} if each point $%
x\in U$ has a neighborhood $W\subset U,$ such that the restriction $%
f|W:W\rightarrow Y$ is (second order) Lipschitz.
\end{itemize}

\begin{Theorem}
\label{T1} Let $X$ be a Banach space and $Y$ be a Banach space with modulus
of convexity of power $2.$ A homeomorphism $f:U\rightarrow V$ between two
open subsets $U\subset X,$ $V\subset Y$ of a Banach space $X$ is locally
convex if

\begin{itemize}
\item the function $f:U\to V$ is locally second order Lipschitz;

\item the function $f^{-1}:V\to U$ is locally Lipschitz.
\end{itemize}
\end{Theorem}

\begin{proof}
Fix any point $x_{0}\in U.$ Given a neighborhood $O(x_{0})\subset U$ of $%
x_{0}\in U\ $ we should construct a convex neighborhood $U(x_{0})\subset
O(x_{0})$ with convex image $f(U(x_{0}))\subset V.$ We lose no generality
assuming that $y_{0}=f(x_{0})=0.$ Since $f$ is locally second order
Lipschitz, the point $x_{0}\in U\ $ has a neighborhood $W\subset O(x_{0})$
such $\omega _{2}(f|W;t)\leq Lt^{2}$ for some real number $L$ and all
positive $t\leq 1.$ Moreover, since the homeomorphism $f^{-1}:V\rightarrow U$
is locally Lipschitz at the point $y_{0}=f(x_{0}),$ we can assume that $%
\omega _{1}(f^{-1}|f(W);t)\leq Lt$ for all $t\geq 0.$ We can, in addition,
also assume that $\max \{\mathrm{diam}(W),\mathrm{diam}f(W)\}\leq 1\}.$

By our assumption, the norm of the Banach space $X$ \ has modulus of
convexity of power type $2\ $and, owing to the relationship \ (\ref{3.3}) of
Lemma \ref{Lm_2}, there is a positive constant $C<1,$ such that
\begin{equation*}
\ \left\Vert \frac{x+y}{2}\right\Vert \leq 1-C\Vert x-y\Vert ^{2}
\end{equation*}%
for any points $x,y\in X$ with $\max \{\Vert x\Vert ,\Vert y\Vert \}\leq 1.$
Take any positive $\varepsilon <1$ such that

\begin{itemize}
\item $\frac{1}{4}L^{2}\leq \frac{C}{\varepsilon };$

\item $B_{\varepsilon }(z_{0})=\{z\in Y:\Vert z-z_{0}\Vert $ $\leq
\varepsilon \}\subset f(W);$
\end{itemize}

The choice of $\varepsilon >0$ guarantees that the preimage $A_{\varepsilon
}=f^{-1}(B_{\varepsilon }(z_{0}))$ of the $\varepsilon $-ball $%
B_{\varepsilon }(z_{0})=\{z\in Y:\Vert z_{0}-z\Vert $ $\leq \varepsilon
\}=\{z\in Y:\Vert z\Vert $ $\leq \varepsilon \}$ centered at $z_{0}=0$ lies
in the neighborhood $W$ on the point $x_{0}\in U.$ The proof of the theorem
will be complete as soon as we check that the closed neighborhood $%
A_{\varepsilon }=f^{-1}(B_{\varepsilon }(z_{0}))\subset W$ of $x_{0}\in U$
is convex. It suffices to check that for any points $x,y\in A_{\varepsilon }$
the point $\bar{x}=\frac{x+y}{2}$ belongs to $A_{\varepsilon }\subset W,$
which happens if and only if its image $f(\bar{x})\in Y$ belongs to the ball
$B_{\varepsilon }(z_{0}).$ Let $h=x-\bar{x}\in X;$ the choice of $%
f:U\rightarrow V$ guarantees that
\begin{equation*}
\begin{array}{c}
\Vert x-y\Vert =\Vert f^{-1}(f(x))-f^{-1}(f(y))\Vert \leq \\
\\
\leq \omega _{1}(f^{-1}|f(W);\Vert f(x)-f(y)\Vert )\leq L\Vert f(x)-f(y)\Vert%
\end{array}%
\end{equation*}%
and
\begin{equation*}
\begin{array}{c}
\left\Vert \frac{f(x)+f(y)}{2}-f(z)\right\Vert =\frac{1}{2}\Vert
f(z+h)-2f(z)+f(x-h)\Vert \leq \\
\\
\leq L\Vert h\Vert ^{2}\leq \omega _{2}(f|W,\Vert h\Vert )\leq L\Vert h\Vert
^{2}=\frac{1}{4}L\Vert x-y\Vert ^{2}\leq \\
\\
\leq \tfrac{1}{4}L^{2}\Vert f(x)-f(y)\Vert ^{2}\leq \frac{C}{\varepsilon }%
\Vert f(x)-f(y)\Vert ^{2}.%
\end{array}%
\end{equation*}%
Since $f(x),f(y)\in B_{\varepsilon }(z_{0}),$ we get $\left\Vert \frac{%
f(x)+f(y)}{2\varepsilon }\right\Vert \ \leq 1-\frac{C}{\varepsilon ^{2}}%
\Vert f(y)-f(x)\Vert ^{2}$ and
\begin{equation*}
\begin{array}{c}
\Vert f(\bar{x})\Vert \text{ }\leq \ \left\Vert \frac{f(x)+f(y)}{2}%
\right\Vert \ +\left\Vert \frac{f(x)+f(y)}{2}-f(\bar{x})\right\Vert \leq \\
\\
\left( \varepsilon -\frac{C}{\varepsilon }\Vert f(y)-f(x)\Vert ^{2}\right) +%
\frac{C}{\varepsilon }\Vert f(y)-f(x)\Vert ^{2}=\varepsilon ,%
\end{array}%
\end{equation*}%
which means that $f(\bar{x})\in B_{\varepsilon }(z_{0})$ and $\bar{x}\in
A_{\varepsilon }\subset f^{-1}(B_{\varepsilon }(z_{0}))\subset W.$
\end{proof}

We say, following \cite{MPVZ}, \ that a function $f:U\rightarrow Y$ defined
on an open subset $U\subset X$ of a Banach space $X$ with values in a Banach
space $Y$ is \emph{G\^{a}teaux differentiable} at a point $x_{0}\in U$ if
there is a linear operator $f_{x_{0}}^{\prime }:X\rightarrow Y$ (called the
\emph{G\^{a}teaux derivative} of $\ f:U\rightarrow Y$ \ at $x_{0}\in U$)
such that for every $h\in X\ $
\begin{equation*}
\lim_{t\rightarrow 0}\frac{f(x_{0}+th)-f(x_{0})}{t}=f_{x_{0}}^{\prime }(h).
\end{equation*}%
By $L(X,Y)$ we denote the Banach space of all bounder linear operators $%
T:X\rightarrow Y$ from $X$ to $Y,$ endowed with the standard operator norm $%
\Vert T\Vert =\sup_{\Vert x\Vert \leq 1}\Vert T(x)\Vert .$ The following
proposition holds.

\begin{Proposition}
Let $X,Y$ be Banach spaces and $U\subset X$ be an open subset. A function $%
f:U\rightarrow Y$ is second order Lipschitz if it is G\^{a}teaux
differentiable at each point of $U$ and the derivative map $f^{\prime
}:U\rightarrow L(X,Y),$ $f^{\prime }:x\mapsto f_{x}^{\prime },$ $x\in X,$ is
Lipschitz.
\end{Proposition}

\begin{proof}
Since the derivative map $f^{\prime }:U\rightarrow L(X,Y)$ is Lipschitz,
there is a constant $L$ such that
\begin{equation*}
\Vert f_{x}^{\prime }-f_{y}^{\prime }\Vert \leq L\Vert x-y\Vert
\end{equation*}%
for each $x,y\in U.$ The second order Lipschitz property of the map $f$ will
follow as soon as we check that
\begin{equation*}
\Vert f(x+h)-2f(x)+f(x-h)\Vert \leq L\Vert h\Vert ^{2}
\end{equation*}%
for each $x\in U$ and $h\in X$ with $[x-h,x+h]\subset U$. The G\^{a}teaux
differentiability of the function $f:U\rightarrow Y$ implies the
differentiability of the function
\begin{equation*}
g:[0,1]\rightarrow Y,\;\;g(t)=f(x+th)-2f(x)+f(x-th).
\end{equation*}%
Moreover,
\begin{equation*}
g^{\prime }(t)=f_{x+th}^{\prime }(h)-f_{x-th}^{\prime }(h)
\end{equation*}%
and hence
\begin{equation*}
\Vert g^{\prime }(t)\Vert =\Vert f_{x+th}^{\prime }-f_{x-th}^{\prime }\Vert
\cdot \Vert h\Vert \leq L\Vert 2th\Vert \cdot \Vert h\Vert =2Lt\Vert h\Vert
^{2}.
\end{equation*}%
Then
\begin{equation*}
\begin{array}{c}
\Vert f(x+h)-2f(x)+f(x-h)\Vert =\Vert g(1)-g(0)\Vert = \\
\\
=||\int_{0}^{1}g^{\prime }(t)dt||\text{ }\leq \int_{0}^{1}\Vert g^{\prime
}(t)||dt=2L\ \Vert h\Vert ^{2}\int_{0}^{1}tdt=L\ \Vert h\Vert ^{2}.%
\end{array}%
\end{equation*}
\end{proof}

This proposition combined with Theorem~\ref{T1} implies

\begin{Corollary}
\label{C1} Let $X$ be a Banach space and $Y$ be a Banach space with modulus
of convexity of power 2. A homeomorphism $f:U\rightarrow V$ between two open
subsets $U\subset X,$ $V\subset Y$ of a Banach space $X$ is locally convex if

\begin{itemize}
\item the function $f:U\rightarrow V$ is G\^{a}teaux differentiable at each
point of $U;$

\item the derivative map $f^{\prime }:U\to L(X,Y)$ is locally Lipschitz;

\item the function $f^{-1}:V\to U$ is locally Lipschitz.
\end{itemize}
\end{Corollary}

The statement of Theorem \ref{T1} is eventually only sufficient. We also at
present do not know if the requirement on the convexity modulus of $Y$ is
essential in Theorem~\ref{T1} and Corollary~\ref{C1}.

\textit{Remark}. Similar to the Hilbert space case, we \ need to mention
here that the integral expressions, considered \ above with respect to the
parameter $t\in \lbrack 0,1],$ are well defined in the Banach space $Y,$ as
the\ related mapping $f_{x(\circ )}^{\prime }(\bar{x}-x_{0}):[0,1]%
\rightarrow Y,$ being continuous and of bounded variation, is \ a priori
Riemann-Birkhoff type integrable \ \ \cite%
{BaMu,BaPo,Birk,CaCr,DiUh,KaSh,Jake}.

\begin{Problem}
\label{prob1} Assume that $X$ is a Banach space such that any locally second
order Lipschitz homeomorphism $f:X\to X$ with locally Lipschitz inverse $%
f^{-1}:X\to X$ is locally convex. Is $X$ (super)reflexive?
\end{Problem}

\begin{Problem}
\label{prob2} For every $n\in\mathbb{N}$ let $F_n$ be the set of all
functions $f_n:l_\infty(n)\to[0,1]$ on the $n$-dimensional Banach space $%
l_\infty(n)=(\mathbb{R}^n,\|\cdot\|_\infty)$ such that

\begin{itemize}
\item $f_{n}^{-1}(0,1]\subset (-1,1)^{n};$

\item $\omega (f;t)\leq t$ for all $t\geq 0;$

\item $\omega _{2}(f;t)\leq t^{2}$ for all $t\geq 0.$
\end{itemize}

Let also $\varepsilon _{n}=\sup \{\Vert f_{n}\Vert _{\infty }:f_{n}\in
F_{n}\}.$ Is $\lim_{n\rightarrow \infty }(1+\varepsilon _{n})^{n}=\infty ?$
\end{Problem}

If \ the Problem~\ref{prob2} has an affirmative answer, then Problem~\ref%
{prob1} has negative answer. Namely, on the reflexive Banach space $%
X=(\sum_{n=1}^{\infty }l_{\infty }(n))_{l_{2}}$ there is a homeomorphism $%
f:X\rightarrow X$ which is not locally convex but $f:X\rightarrow X$ \ is
second order Lipschitz and $f^{-1}:X\rightarrow X$ is locally Lipschitz.

\section{Conclusion}

The paper analyses of smooth nonlinear mappings for Hilbert and Banach spaces that carry small balls to convex sets, provided that the radius of the balls is sufficiently small. The main goal is to establish new and mild sufficient conditions for a nonlinear mapping to be locally convex. The analysis involves both Hilbert and Banach spaces, and some of the results are found to be interesting and novel even for finite-dimensional problems.

We specifically address a suitably reformulated local convexity problem for both Hilbert and Banach spaces. The local convexity property holds for differentiable mappings of Hilbert spaces if the Frech\'{e}t derivative is Lipschitzian in a closed ball with certain properties. We provide an improved estimation for the radius of the ball, ensuring that its image is convex. This result is established with arguments different from previous works, leading to a more refined analysis.

For Banach spaces, the local convexity problem is more intricate and requires subtle techniques. We analyze locally convex functions between Banach spaces, considering the interplay between the modulus of convexity of a Banach space and the modulus of smoothness of a function. This generalization allows for a deeper understanding of the local convexity property in the Banach space setting. Some of the presented results are novel even in the finite-dimensional case.

Moreover, we formulate open problems related to the local convexity property for nonlinear mappings of Banach spaces, highlighting areas where further research is needed.

\section{Acknowledgements}

The authors are grateful to A. Augustynowicz, D. Blackmore, L. G\'{o}%
rniewicz and A. Plichko for fruitful discussions and remarks. They are
especially indebted to T. Banakh for generous and invaluable help in
treating the Banach space case and mentioning the references related with
the topic studied in the article.

\subsection*{Contributions}

Conceptualization, Y.P., P.P.; methodology M.V., M.G.; validation, Y.P., P.P.; investigation, Y.P., P.P.; writing original draft preparation, M.V., M.G.; writing-review and editing, M.V., M.G.; project administration, Y.P., P.P.; funding acquisition, Y.P., P.P. All authors have read and agreed to the published version of the manuscript.

\subsection*{Funding}

The results were obtained as part of the work under a grant from the Ministry of Education and Science of Ukraine (project number 0123U101691).

\subsection*{Data Availability}

Not applicable.

\subsection*{Conflict of interest}

Authors report no potential conflict of interest.

\begin{adjustwidth}{-\extralength}{0cm}

\reftitle{References}

\PublishersNote{}
\end{adjustwidth}

\begin{thebibliography}{99}
\bibitem{AlTiFo} Alexeev, V.M.; Tikhomirov, V.M.; Fomin S.V. . \textit{Optimal Control};
Consultants Bureau: New York, USA, 1987.

\bibitem{Aug} Augustynowicz, A.; Dzedzej, Z.; Gelman, B.D. The
solution set to BVP for some functional-differensial inclusions. \textit{Set-Valued Analysis}. {\bf 1998},\textit{6}, 257--263. https://doi.org/10.1023/A:1008618606813

\bibitem{BaBaPlPr} Banakh, I.; Banakh, T.; Plichko, A\.; Prykarpatsky, A. On local convexity of nonlinear mappings between Banach spaces. \textit{Cent. Eur. J.
Math.} {\bf 2012}, \textit{10(6)}, 2264--2271.  https://doi.org/10.2478/s11533-012-0101-z

\bibitem{BaSa} Baccari, A.; Samet, B. An extension of Polyak's theorem in a Hilbert space. \textit{J. Optim. Theory Appl.} {\bf 2009}, \textit{140(3)}, 409--418. https://doi.org/10.1007/s10957-008-9457-4

\bibitem{BoEmKo} Bobylev, N.A.; Emelyanov, S.V.; Korovin, S.K. Convexity of images of convex sets under smooth maps. Nonlinear Dynamics and Control. \textit{Comput. Math. Model.} {\bf 2009}, \textit{15(3)}, 213--222. https://doi.org/10.1023/B:COMI.0000035819.33749.a7

\bibitem{Che} Chernukha, O.; Chuchvara, A.; Bilushchak, Y.; Pukach, P.; Kryvinska, N.
Mathematical Modelling of Diffusion Flows in Two-Phase Stratified Bodies with Randomly Disposed Layers of Stochastically Set Thickness. \textit{Mathematics}. {\bf 2022}, \textit{10(19)}, 3650. https://doi.org/10.3390/math10193650

\bibitem{Dyma} Dymarsky, A. Convesity of a small ball under quadratic map. \textit{Linear Algebra Appl.} {\bf 2016}, \textit{488}, 109--123. \\  https://doi.org/10.1016/j.laa.2015.09.005

\bibitem{BaPo} Balcerzak, M.; Potyrala, M. Convergence theorems for the
Birkhoff integral. \textit{Czechoslovak Mathematical Journal}. {\bf 2008}, \textit{58(4)}, 1207--1219. Available at:  http://dml.cz/dmlcz/140451 http://dml.cz/dmlcz/140451

\bibitem{BaMu} Balcerzak, M.; Musial, K. A convergence theorem for the Birkhoff integral. \textit{Funct. Approx. Comment. Math.} {\bf 2014}, \textit{50(1)}, 161--168. https://doi.org/10.7169/facm/2014.50.1.5

\bibitem{BP} Blackmore, D.; Prykarpatsky, A.K. A solution set analysis of a nonlinear operator equation using a Leray Schauder type fixed point approach. \textit{Topology}. {\bf 2009}, \textit{48}, 182--185. https://doi.org/10.1016/j.top.2009.11.017

\bibitem{Birk} Birkhoff, G. Integration of functions with values in a Banach
space. \textit{Trans. Amer. Math. Soc.} {\bf 1935}, \textit{38}, 357--378. https://doi.org/10.2307/1989687

\bibitem{BGHV} Borwein, J.; Guirao, A.; H\'{a}jek; P., ~Vanderwerff, J. Uniformly convex functions on Banach spaces. \textit{Proc. Amer. Math. Soc.} {\bf 2009}, \textit{137(3)}, 1081--1091. Available at:
https://www.ams.org/journals/proc/2009-137-03/S0002-9939-08-09630-5/S0002-9939-08-09630-5.pdf

\bibitem{CaCr} Candeloro, D.; Croitoru, A.; Gavrilute, A.; Iosif, A.; Sambucini, A.R. Properties of the Riemann-Lebesgue integrability in the
non-additive case. \textit{arXiv:1905.03993v1 [math.FA]}. {\bf 2019}. https://doi.org/10.48550/arXiv.1905.03993

\bibitem{DiUh} Diestel, J.; Uhl (Jr.), J. J. \textit{Vector measures}; Math. Surveys Monogr., 15, Amer. Math. Soc.: Providence, RI, 1977. \\ https://doi.org/10.1090/surv/015

\bibitem{FaHa} Fabian, M.F.; Habala, P.; H\'{a}jek, P.; Montesinos, S.V.; Pelant, J.;  Zizler, V. \textit{Functional Analysis and Infinite-Dimensional
Geometry}; Springer-Verrlag: New York, USA, 2001. https://doi.org/10.1007/978-1-4757-3480-5

\bibitem{Frid-1} Frid, H. Nonlinear maps of convex sets in Hilbert spaces
with application to kinetic equations. \textit{Bull. Braz. Math. Soc., New Series.} {\bf 2008}, \textit{39(3)}, 315--340. https://doi.org/10.1007/s00574-008-0008-2

\bibitem{Frid-2} Frid, H. \textit{Maps of convex sets in Hilbert spaces}; Universit%
\"{a}t G\"{o}ttingen, 2008. \\ Available at:  https://webdoc.sub.gwdg.de/ebook/serien/e/IMPA-A/414.pdf

\bibitem{Ge1} Goebel, K. \textit{Zagadnienia metrycznej teorii punkt\'{o}wstalych}; Wydawnictwo Uniwersytetu Marii-Curie Sklodowskiej: Lublin, 1999. (in Polish)

\bibitem{Ge2} Goebel, K. \textit{Twierdzenia o punktach stalych. Wyklady}; Wydawnictwo Uniwersytetu Marii-Curie Sklodowskiej: Lublin, 2005. (in Polish)

\bibitem{Go} G\'{o}rniewicz, L. \textit{Topological fixed point theory of
multi-valued mappings}; Kluwer: Dordrecht, 1999. https://doi.org/10.1007/978-94-015-9195-9

\bibitem{Greg} Gregus. M. A fixed point theorem in Banach spaces. \textit{Boll. Un. Math. Ital.} {\bf 1980}, \textit{5(17)}, 193--198.

\bibitem{GH} Guirao, A.J.; H\'{a}jek, P. On the moduli of convexity. \textit{Proc. Amer. Math. Soc.} {\bf 2007}, \textit{135(10)}, 3233--3240. Available at: http://www.jstor.org/stable/20534944

\bibitem{HMZ} H\'{a}jek, P.; Montesinos, V.; ~Zizler, V. Geometry and Gateaux smoothness in separable Banach spaces.  \textit{Operators and Matrices}. {\bf 2012}, \textit{6(2)}, 201--232.  https://doi.org/10.7153/oam\%2D06\%2D15

\bibitem{Ho} H\"{o}rmander, L. Sur la fonction d'applui des ensembles convexes dans une espace localement convexe. \textit{Arkiv Math.} {\bf 1955}, \textit{3(2)}, 180--186. https://doi.org/10.1007/BF02589354

\bibitem{IoTi} Ioffe, A.D., Tikhomirov, V.M. \textit{Theory of extremal problems}; North-Holland Publ. Co.: Amsetdam-New York, 1979.

\bibitem{Ivan} Ivanov, G.E. Nonlinear images of sets. I. Strong and weak
convexity. \textit{J. Convex Anal.} {\bf 2020}, \textit{27(1)}, 363--382.

\bibitem{Jake} Jaker, Ali Sk. Riemann and Riemann type integration in
Banach spaces. \textit{Real Analysis Exchange}, {\bf 2013/14}, \textit{39(2)}, 403--440.

\bibitem{Ledy} Ledyaev, Yu.S. Criteria for the convexity of closeds sets in
Banach spaces. \textit{Proc. Steklov Math. Inst.} {\bf 2019}, \textit{304(1)}, 190--204. https://doi.org/10.1134/S0081543819010139

\bibitem{Matv} Matviychuk, Y.; Kryvinska, N.; Shakhovska, N., Poniszewska-Maranda, A. New
principles of finding and removing elements of mathematical model for reducing computational and time complexity.  \textit{International Journal of Grid and Utility Computing}. {\bf 2023}, \textit{14(4)}, 400--410.   https://doi.org/10.1504/ijguc.2023.132625

\bibitem{Mord} Mordukhovich, B.S. \textit{Variatioonal analysis and generalized
differentiation.I. Basic thepry}; Springer-Verlag: Berlin, Germany, 2006. https://doi.org/10.1007/3-540-31247-1

\bibitem{Naz} Nazarkevych, M.; Kryvinska, N.; Voznyi, Y. Applying ateb-gabor filters to
biometric imaging problems. \textit{Symmetry}. {\bf 2021},  \textit{13(4)}, 717  https://doi.org/10.3390/sym13040717

\bibitem{KaSh} Kadets, V.M.; Shumyatskiy, B.; Shvidkoy, R.; Tseytlin, L.M.; Zheltukhin, K. Some remarks on vector-valued integration. \textit{Mat. Fiz. Anal.
Geom.} {\bf 2002}, \textit{9}, 48--65. https://hdl.handle.net/11511/76330

\bibitem{KZ} Krasnoselsky, M.A.; Zabreyko, P.P. \textit{Geometric methods
of nonolinear analysis}; Nauka: Moscow, USSR, 1975. (in Russian)

\bibitem{Li} Linke, Y.E. Application of Michael's theorem and its
converse to sublinear operators. \textit{Mathematical Notes}, {\bf 1992}, \textit{52(1)}, 680--686. https://doi.org/10.1007/BF01247650

\bibitem{MPVZ} McLaughlin, D.; Poliquin, R.; Vanderwerff, J.; Zizler, V. Second-order Gateaux differentiable bump functions and approximations in
Banach spaces. \textit{Canad. J. Math.}  {\bf 1993}, \textit{45(3)}, 612--625. https://doi.org/10.4153/CJM-1993-032-9

\bibitem{Nir} Nirenberg, L. \textit{Topics in Nonlinear Functional Analysis}; AMS Publisher: New York, USA, 1974.

\bibitem{PP} Petunin. Y.I.; Plichko, A.I. \textit{Characteristics Theory for
Subspaces and Its Applications.}; Vyshcha Shkola: Kyiv, Ukraine, 1980. (in Russian)

\bibitem{Phel} Phelps, R.R. Convex Functions, Monotone Operators and
Differentiability. \textit{Lecture Notes in Mathematics}. {\bf 1993}, \textit{1364}. https://doi.org/10.1007/978-3-540-46077-0

\bibitem{Po} Polyak, B.T. Convexity of nonlinear image of a small ball with
applications to optimization. \textit{Set-Valued Analysis}.  {\bf 2001}, \textit{9}, 159--168. https://doi.org/10.1023/A:1011287523150

\bibitem{PBPP} Prykarpatska, N.K.; Blackmore, D.L.; Prykarpatsky, A.K.; Pytel-Kudela, M. On the inf-type extremality
solutions to Hamilton--Jacobi equations and some generalizations. \textit{
Miskolc Math. Notes}. {\bf 2003}, \textit{4(2)}, 153--176. https://doi.org/10.18514/MMN.2003.66

\bibitem{Pr1} Prykarpatsky, A.K. An infinite dimensional
Borsuk-Ulam type generalization of the Leray-Schauder fixed point theorem
and some applications. \textit{Ukrainian Mathematical Journal}. {\bf 2008}, \textit{60(1)}, 114--120. https://doi.org/10.1007/s11253-008-0046-3

\bibitem{Pr} Prykarpatsky, A.K. A Borsuk Ulam type generalization of
the Leray-Schauder fixed point theorem.  \textit{arXiv:0902.4416}. {\bf 2009}.  https://doi.org/10.48550/arXiv.0902.4416

\bibitem{Pryk} Prykarpatsky, A.K. \textit{A Borsuk-Ulam type generalization of the
Leray-Schauder type fixed point theorem}; United Nations Educational,
Scientific and Cultural Organization and International Atomic Energy Agency, 2007. Available at: \\
https://inis.iaea.org/collection/NCLCollectionStore/ Public/39/016/39016745.pdf?r=1

\bibitem{Reis} Reissig, G. Convexity of the reachable sets of control
systems. \textit{Autom. Remote Control}. {\bf 2007}, \textit{68(9)}, 1527--1543.\\ https://doi.org/10.1134/S000511790709007X

\bibitem{SPS} Samoilenko, A.M.; Prykarpats'kyi, A.K.;  Samoilenko, V.H.
Lyapunov--Schmidt approach to studying homoclinic splitting in
weakly perturbed Lagrangian and Hamiltonian systems. \textit{Ukr. Mat. Zh.}.  {\bf 2003}, \textit{55(1)}, 82--92. https://doi.org/10.1023/A:1025072619144

\bibitem{Sch} Schwartz, J.T. \textit{Nonlinear functional analysis}; Gordon
and Breach Science Publisher: New York, USA, 1969.

\bibitem{Vovk} Vovk, M.I.; Pukach, P.Ya.; Dilnyi, V.M.; Prykarpatski, A.K. Hilbert Space Decomposition Properties of Complex Functions and Their Applications. \textit{Contemporary Mathematics (Singapore)}. {\bf 2023}, \textit{4(4)}, 702--709.  https://doi.org/10.37256/cm.442023386
\end{thebibliography}
\end{document}